%% file: main.tex
\documentclass[letterpaper, 10 pt, conference]{ieeeconf}  % Comment this line out
                                                       
\IEEEoverridecommandlockouts                              % This command is only
\usepackage{amsmath} % assumes amsmath package installed
\input{macro/jared_macro}

\input{macro/data_driven_macro}
\usepackage[colorinlistoftodos]{todonotes}
\usepackage{pmat}

\usepackage{color}
\usepackage{ulem}

\title{\LARGE \bf Data-Driven Structured Robust Control of Linear Systems
}

\author{Jared Miller$^1$,  Jaap Eising$^1$, Florian D\"{o}rfler$^1$, Roy S. Smith$^1$
\thanks{$^1$J. Miller, J. Eising, F. D\"{o}rfler, and R. S. Smith are with the Automatic Control Laboratory (IfA), ETH Z\"{u}rich, Switzerland (\{jarmiller, jeising, doerfler, rsmith\}@control.ee.ethz.ch).}
\thanks{This work is supported in part by NCCR Automation, funded by the Swiss National Science Foundation (grant number 51NF40\_225155) and the SNF/FW Weave Project 200021E\_20397.
}
}

\begin{document}

\maketitle
\thispagestyle{empty}
\pagestyle{empty}

%%%%%%%%%%%%%%%%%%%%%%%%%%%%%%%%%%%%%%%%%%%%%%%%%%%%%%%%%%%%%%

\input{sections/abstract}
\input{sections/introduction}

% \input{sections/contributions}
\input{sections/summary}
\input{sections/preliminaries}
\input{sections/data_driven_structure}
\input{sections/examples}
\input{sections/Conclusion}
\input{sections/acknowledgements}

% \input{sections/appendix}

% %%%%%%%%%%%%%%%%%%%%%%%%%%%%%%%%%%%%%%%%%%%%%%%%%%%%%%%%%%%%%%%%%%%%%%%%%%%%%%%%
% \section{Acknowledgements}

% The authors thank Milan Korda for his discussions about occupation measures and time-varying uncertainty.

% \todo{fix references}

\bibliographystyle{IEEEtran}
\bibliography{references.bib}

\end{document}

%% file: macro/jared_macro.tex
\usepackage[nolist]{acronym}
\usepackage{xcolor}
\usepackage{amsmath}
\usepackage{amssymb}
\usepackage{mathtools}
\usepackage{subcaption}
\usepackage{hyperref}
\hypersetup{colorlinks=true, linkcolor=black}
\usepackage[ruled]{algorithm2e}
\usepackage{mathrsfs}
\definecolor{ccolor}{RGB}{203,96,21}

\newcommand{\R}{\mathbb{R}}

\newcommand{\dc}{\mathcal{D}}

\newcommand{\s}{\mathcal{S}}

\newcommand{\0}{\mathbf{0}}
\newcommand{\1}{\mathbf{1}}
\newcommand{\psd}{\mathbb{S}}

\DeclarePairedDelimiter{\norm}{\lVert}{\rVert}

\DeclarePairedDelimiterX{\inp}[2]{\langle}{\rangle}{#1, #2}

\newtheorem{thm}{Theorem}[section]
\newtheorem{lem}[thm]{Lemma}

\newtheorem{exmp}{Example}[section]
\newtheorem{rmk}{Remark} 
\newtheorem{asmp}{Assumption}

%% file: macro/data_driven_macro.tex
\newcommand{\bw}{\mathbf{W}_-}

\newcommand{\xp}{\mathbf{X}_+}
\newcommand{\xn}{\mathbf{X}_-}

\newcommand{\bu}{\mathbf{U}_-}

\usepackage{harmony}

%% file: sections/abstract.tex
\begin{abstract}
\label{sec:abstract}
Static structured control refers to the task of designing a state-feedback controller such that the control gain satisfies a subspace constraint. Structured control has applications in control of communication-inhibited dynamical systems, such as systems in networked environments. This work performs  $H_2$-suboptimal regulation under a common structured state-feedback controller for a class of data-consistent plants. The certification of $H_2$-performance is attained through a combination of standard $H_2$ LMIs, convex sufficient conditions for structured control,  and a matrix S-lemma for set-membership. The resulting convex optimization problems are linear matrix inequalities whose size scales independently of the number of data samples collected. Data-driven structured $H_2$-regulation control is demonstrated on example systems.
% sufficient conditions for a class of plants compati

\end{abstract}

%% file: sections/introduction.tex
\section{Introduction}
\label{sec:introduction}

When designing control policies for dynamical system, it is often desired that the controller respect a given structural property.
In a state-feedback context, the gain matrix $K$ of the control policy $u = K x$ must respect a subspace constraint $K \in \s$, in which $\s$ encodes the desired structure of the problem. 
Structured control with subspace constraints is generically an NP hard problem \cite{blondel1997np}.
A specific instance of structured control is decentralized control, in which the controller must match a sparsity pattern derived from the information/actuation structure of the overall system \cite{wang1973stabilization, bakule2008decentralized}. 
% , even for the decentralized setting.
Decentralized control can be convex in the case where the subspace $\s$ possesses a property of \textit{quadratic invariance} (QI) \cite{lessard2011qi}, but QI-based convex control schemes are only nonconservative if dynamic feedback laws are permitted (rather than static output feedback). Network decompositions can be performed to reduce the conservatism of decentralized control methods, which would in turn be nonconservative under dynamic output feedback and QI \cite{furieri2018robust}. The work in \cite{furieri2019separable} performs iterative merging to create separable quadratic Lyapunov functions in the distributed (static state feedback) control setting. The work in \cite{ferrante2019design, ferrante2020lmi} offer convex but conservative formulations for state-feedback structured control in the stabilization and $H_\infty$-norm minimization settings.

% Sparsity-based $\s$ patterns are a specific instance of sign-constrained control (ensuring that $K$ satisfies a pattern of signs). Sign-based structure can be imposed through affine sign constraints when searching over polytopic Lyapunov functions (e.g. superstability) \cite{polyak2002superstable}. A specific instance of this polytopic search includes optimizing over (dual) linear copositive Lyapunov functions for positive (or cone-invariant) dynamical systems \cite{rantzer2018tutorial}.

This work focuses on finding a worst-case $H_2$-regulating structured control for all systems that are compatible with observed data using the framework of set-membership data-driven control. The methodology of data-driven control treats the data itself as a model of the system, and then synthesizes control parameters based on the data without first passing through a system identification routine \cite{HOU20133, hou2017datasurvey}. The idea behind such approaches is the fact that finding a control law is often easier than describing a suitable model for a plant. One approach is the set-membership style of data-driven control, which requires describing:
\begin{enumerate}
    \item The set of plants regulated by a given control policy;
    \item The set of plants consistent with the observed data under an a-priori-known noise bound;
    \item A certificate of set containment such that every plant consistent with the data is also regulated.
\end{enumerate}

This paper will describe the set of plants consistent with the data by a matrix ellipsoid as defined by a \ac{QMI} constraint \cite{van2023quadratic}. The work in \cite{van2020noisy} demonstrates how to derive guarantees for stabilization, worst-case $H_2$ suboptimal control, and worst-case $H_{\infty}$ suboptimal control over these \ac{QMI}-defined sets through the use of a matrix S-Lemma \cite{yakubovich1997s}. These \ac{QMI}-based results have also been used for control of networked systems under block-diagonal partitioning in \cite{eising2022informativity,    wang2023event} and for pole placement in LMI-defined regions ($D$-stability) \cite{bisoffi2022learning}. Other instances of set-membership data-driven control include the formation of polytopic Lyapunov functions (e.g. superstability) \cite{cheng2015, miller2023positive} with certifications arising from an Extended Farkas Lemma \cite{hennet1989farkas}.

Methodologies for data-driven control that require a model reference include virtual reference feedback tuning \cite{campi2002virtual}, \cite{bazanella2011data}, iterative feedback tuning \cite{hjalmarsson1998iterative}, and correlation-based tuning \cite{karimi2004iterative}. A particularly popular technique for data-driven control is the behavioral approach based on Willem's Fundamental Lemma \cite{willems2005note}, in which the space of all feasible trajectories of a linear system can be expressed as a linear combination of observed data (in Hankel matrix form) if a persistency of excitation rank condition is obeyed. This Fundamental Lemma has been employed for model predictive control
\cite{waarde2020informativity, depersis2020formulas, coulson2019data, berberich2020mpc}, stochastic control \cite{pan2022stochastic}, and closed-loop control 
\cite{DINKLA20231388}.
% \cite{DINKLA20231388, moffat2024transient}.
The work of \cite{celi2023data} performs sparse (structured) control design of a feedback gain $K$ in the behavioral setting under the absence of noise.

The contributions of this paper are:
\begin{enumerate}
    \item A centralized LMI formulation for data-driven $H_2$ regulation using structured control based on the sufficient-convex constraint of \cite{ferrante2019design};
    \item An accounting of computational complexity for the presented conditions;
    \item Demonstrations of $H_2$-suboptimal structured control on example systems.
\end{enumerate}

% Behavioral theory originated from the work of Jan C. Willems, who described a system as the collection of possible generated temporal trajectories. The behavior of \ac{LTI} systems is completely specified by the column span of the input-output Hankel matrix if a persistency of excitation rank condition is satisfied (Willem's Fundamental Lemma) \cite{willems2005note}. The Fundamental Lemma has been used for predictive control, stabilization, and worst-case-optimal control \cite{waarde2020informativity, depersis2020formulas, coulson2019data, berberich2020mpc}. 

% $H_2$-suboptimal structured control will be performed by adapting methods used in the set-membership data-driven control literature for unstructured $H_2$-suboptimal robust control. 
% \urg{Survey of set-membership-style DDC and other methods.}
% The set of data-consistent systems will be represented as a matrix ellipsoid (quadratic matrix inequality region), and the certificate of set-containment will be enforced through a Matrix S-Lemma \cite{van2020noisy}. 

% \urg{more citations of data-driven control. Include Jaap's networked results in \cite{eising2022informativity}, which was used in the event-triggered framework \cite{wang2023event}.}

%% file: sections/summary.tex
This paper has the following structure: 
Section \ref{sec:preliminaries} reviews preliminaries such as notation and (nominal) structured control; Section \ref{sec:data_driven} describes consistency of structured-stabilizing-plants and data-consistent plants, and develops an \ac{LMI} for data-driven certifiable $H_2$-suboptimal structured control;
% Section \ref{sec:h2} modifies the \ac{LMI} to demonstrate the capability of $H_2$-suboptimal structured control \urg{if space allows}. 
Section \ref{sec:examples} performs demonstrations on example systems; and Section \ref{sec:conclusion} concludes the paper.
% \urg{Fill in the paper structure}
% Section \ref{sec:preliminaries} will review preliminaries such as notation, notions of stability for linear systems, and \ac{SOS} proofs of polynomial nonnegativity. Section \ref{sec:full_method} will present 
% The paper is concluded in Section \ref{sec:conclusion}.

%% file: sections/preliminaries.tex
%\section{Preliminaries}

% \subsection{Acronyms/Initialisms}
\input{sections/acronym}

\subsection*{Notation}

The set of $m \times n$ real-valued matrices is $\R^{m \times n}$. The transpose of a matrix $M \in \R^{m \times n}$ is $M^\top \in \R^{n \times m}$. The set of $n \times n$ square symmetric matrices is denoted $\psd^n$.
%(matrices with $M = M^\top$)
If $P \in\psd^n$, we write $P\succeq 0$, and $P \succ 0$ if P is positive semidefinite or positive definite, respectively. The sets of positive semidefinite and positive definite matrices are denoted $\psd^n_+$ and $\psd^n_{++}$ respectively. 
%The doubled-Hermitian operator of a square matrix $M \in \R^{n \times n}$ is $\text{He}(M) = M + M^\top$. The notation $\star$ will refer to canonical transposed components in a block matrix such as \todo{Near the end, check whether we can avoid introducing He and $\star$. Less abbr. makes things more legible.} 
%\begin{align}
%    \begin{bmatrix}
%        P & Q \\
%        Q^\top &  I
%    \end{bmatrix} &\rightarrow \begin{bmatrix}
%        P & Q \\
%        \star &  I
%    \end{bmatrix}
%\intertext{The $\star$ notation will also be used to abbreviate canonical transposed elements in quadratic expressions as in}
%    Z M Z^\top &\rightarrow Z M [\star]^\top.
%\end{align}
% \todo{Add kronecker product.}
The Kronecker product between matrices $A$ and $B$ is $A \otimes B$.
The symbols $I_n$, $\0_{n \times m}$ and $\1_{n \times m}$ will refer to the identity, zeros, and ones matrix respectively. The dimensions will be omitted if there is no ambiguity. Given a transfer matrix $G(z)$, we denote its $H_2$-norm by $\norm{G(z)}_{H_2}$. 
% The transition of a discrete time system for a state $x$ will be denoted as $x_+$ (such as $x_+ = A x$ for a linear system $x_{t+1} = A x_t$).

\section{Preliminaries}

\label{sec:preliminaries}

%\subsection{LMIs for $H_2$ Control Synthesis}
% This subsection will summarize the work of \cite{ferrante2019design} for stabilizing structure control. The structured control methodology has been extended in \cite{ferrante2020lmi} to $H_\infty$-minimizing structured control (with the requirement that a line search must be performed over a single free parameter).

% For simplicity, this section will focus on the discrete-time setting

This work considers discrete-time linear systems with state $x \in \R^n$, input $u \in \R^m$, exogenous input $\xi \in \R^z$, and measured output $y \in \R^q$. The next state and output at each discrete time index $t$ are 
% \todo{I saw some conflicting notes re: $x_+, x^+$ etc. Why not go with the form of \eqref{eq:stab_W} immediately? I don't think it'll cost us any space. }
\begin{align}\label{eq:sys_h2} x(t+1) &= A x(t) + B u(t) + E \xi(t)  \\
y(t) &= Cx(t) + D u(t). \nonumber \end{align} 
Given a static state-feedback controller $u=Kx$, consider the closed-loop system 
\begin{align}  x(t+1) &= (A+BK)x(t)+E\xi(t)  \\y(t)&= (C+DK)x(t), \end{align}
and denote the resulting closed-loop transfer function from $\xi$ to $y$ by $G_K(z)$. The classical suboptimal $H_2$ control design problem (see e.g. \cite[Chapter 4] {zhou1998essentials}) reads as follows: Given a \textit{performance metric} $\gamma\geq 0$, design (if it exists) $K \in \R^{n \times m}$ such that $\norm{G_K(z)}_{H_2}\leq\gamma$. In this case, we say $K$ is a $\gamma$-suboptimal $H_2$ controller for \eqref{eq:sys_h2}. The $H_2$ synthesis procedure can be readily cast as an \ac{LMI} feasibility problem (see e.g. \cite{boyd1994linear}). Here, we follow the approach of \cite{de2002extended}, which introduces an additional variable $R$, which will prove useful in the remainder of the paper. 

\begin{lem}[Theorem 1 of \cite{de2002extended}]
    \textit{The feedback gain $K$ is a $\gamma$-suboptimal $H_2$ controller for \eqref{eq:sys_h2} if there exists matrices $(P, Q, R)$ such that the following holds: }
    % \todo{(1) note that this is not an LMI, (2) we also need P>0 right? We do technically get it from (a) I guess. (3) isn't this an if and only if statement?}
    \begin{subequations}   
    \label{eq:h2_basic}
    \begin{align}
        %\find_{P, Q, R, K} 
        &\begin{bmatrix}
            P & (A + B K) R  & E\\
            R^\top (A + B K)^\top & R + R^\top - P & 0 \\
            E^\top & 0 & I
        \end{bmatrix} \succ 0 \label{eq:eet_basic}\\
        & \begin{bmatrix}
            Q & (C+DK) R\\ R^\top (C+DK)^\top & R + R^\top - P
        \end{bmatrix} \succ 0 \\
        & \text{Tr}(Q) \leq \gamma^2 \\
        & P \in \psd^n, \ Q \in \psd^q, \ R \in \R^{n \times n}, K\in\R^{m\times n}.
    \end{align}
    \end{subequations}
\end{lem}

% \begin{prop}
%     The following statements are equivalent for a continuous-time system
%         \begin{align}
%         & \exists P \in \psd^n_{++}:  \text{He}((A+BK) P \prec 0)\label{eq:design_k_only}\\
%                 & \exists P \in \psd_{++}^n, R \in \R^{n \times n}: \nonumber \\
%                 & \qquad \qquad\begin{bmatrix}
%             0 & P \\
%             P & 0
%         \end{bmatrix} + \text{He}\left(\begin{bmatrix}
%             A+BK \\ -I
%         \end{bmatrix}\begin{bmatrix}
%             R & R
%         \end{bmatrix}\right) \prec 0. \label{eq:design_kr_c}
%     \end{align}
    
%     and for a discrete-time system (Equation (10) of \cite{pipeleers2009extended}):
%     \begin{align}
%         & \exists P \in \psd^n_{++}: \quad & &\begin{bmatrix}
%             P & (A+BK) P \\
%             \star & P
%         \end{bmatrix} \succ 0 \label{eq:design_k_only_d}\\
%         & \exists P \in \psd_{++}^n, R \in \R^{n \times n}: \quad & &\begin{bmatrix}
%             P & (A+BK) R \\
%             \star & R + R^\top - P
%         \end{bmatrix} \succ 0. \label{eq:design_kr_d}
%     \end{align}
% \end{prop}
Since \eqref{eq:h2_basic} is not linear in all decision variables, it cannot be resolved efficiently. However, we can define a new variable $L = K R$, which leads to a linear problem. 
\begin{lem}[Theorem 5 of \cite{de2002extended}] \label{lem:h2}
% A $\gamma$-suboptimal  $H_2$
\textit{There exist matrices $(P,Q,R,K)$ such that \eqref{eq:h2_basic} holds if and only if the following LMI in $(P, Q, R, L)$ is feasible}
% A controller $u=Kx$ exists regulating the system in \eqref{eq:lin_sys} with $H_2$ norm less than $\gamma$ \eqref{eq:lin_sys} if the following \ac{L}
%     \inf_{P, R, L, Q, \gamma } & \qquad \gamma \\
\begin{subequations}
    \label{eq:h2_design}
\begin{align}
        \quad  &\begin{bmatrix}
            P - E E^\top & A R + B L \\
            (A R + B L)^\top & R + R^\top - P 
        \end{bmatrix} \succ 0 \label{eq:eet_schur}\\
        & \begin{bmatrix}
            Q & CR +D L\\ (CR +D L)^\top & R + R^\top - P
        \end{bmatrix} \succ 0 \\
        & \text{Tr}(Q) \leq \gamma^2 \\
        & P \in \psd^n, \ Q \in \psd^q, \ R \in \R^{n \times n}, \ L \in \R^{m \times n}.
    \end{align}
\end{subequations}
\textit{Moreover, if the LMIs in \eqref{eq:h2_design} are feasible, then the matrix $R$ will be nonsingular and $K = L R^{-1}$ is a $\gamma$-suboptimal $H_2$ controller.}
\end{lem}
In addition to finding $\gamma$-suboptimal $H_2$ controllers for fixed $\gamma$, we are also interested in infimizing the error bound $\gamma$ (or, maximizing the performance). For this, note that the previous matrix inequality is linear in all its decision variables and in $\gamma^2$. Hence, finding $K$ such that $\gamma$ is minimal takes the form of a semidefinite program (SDP) and can be solved with standard methods. 
% To be precise, \acp{LMI} \eqref{eq:eet_basic} and \eqref{eq:eet_schur} are related through the substitution $L = K R $ and by applying a Schur complement.
%An advantage of the formulation in \eqref{eq:h2_design} as compared to \eqref{eq:h2_basic} is that the Lyapunov matrix $P$ can be chosen to be dense even while $K$ and $R$ are sparse. This freedom will be used to synthesize structured controllers.

\subsection{Convex Sufficient Conditions for Structured Control}
For many applications, full state feedback controllers are not desirable or are even impossible to implement. As an example, networked systems or networked state-feedback may have limited communication ability, which in turn imposes a sparsity pattern on the set of feasible controller matrices $K$.  Here, we consider the problem of \textit{structured control design}. That is, we restrict ourselves to feedback gains $K\in\mathcal{S}$, where $\mathcal{S} \subseteq \mathbb{R}^{m\times n}$ is a given subspace.

\begin{exmp}[Multi-agents and sparsity patterns]
\label{exmp:sparsity_explain}
    One common source of subspace constraints in $\s$ are \textit{sparsity constraints}. This can be captured in a prescribed elementwise zero-nonzero structure on $K$. The enforced zeros in certain elements of the $K$ matrix could correspond to lack of information or actuation capacity between output and input channels, such as in the networked or multi-agent setting. An example of a sparsity constraint $\s$ is:
    \begin{equation} \label{eq:sparse} \s = \left\lbrace \begin{bmatrix} \alpha_1 & \alpha_2 & 0 \\ 0 & \alpha_3 & \alpha_4 \end{bmatrix}  \quad \Bigg|  \quad  \alpha_1,\ldots,\alpha_4\in\R \right\rbrace. \end{equation} 
    If the subspace $\s \subseteq \mathbb{R}^{m\times n} $ arises from a sparsity pattern, we will identify it with \textit{sparsity pattern} $\text{sp}(\s)\in \{0,1\}^{m\times n}$, where $0$ denotes an enforced zero and $1$ a free element. For example, for \eqref{eq:sparse} this yields
    \begin{align}
         \text{sp}(\s) = \begin{bmatrix}
            1 & 1 & 0 \\ 0 & 1 & 1
        \end{bmatrix}.
    \end{align}
\end{exmp}

Looking back at Lemma~\ref{lem:h2}, the constraint $K = L R^{-1} \in \s$ is generically \textit{not} convex in $L$ and $R$. A common approach for subspaces $\s$ defined by a sparsity pattern, is to impose that $L \in \s$ and $R$ is diagonal. 
% An advantage of the formulation in 
In contrast, the work in \cite{ferrante2019design} poses a less conservative sufficient condition to ensure $K \in \s$. The approach in \cite{ferrante2019design} generalizes to subspaces $\s$ that do not arise from sparsity patterns, such as the setting of coordinated control in a multi-agent system. To describe the subspace structure, we first require some notation. Given a basis $\{S_\ell\}_{\ell=1}^k$ for $\s$ (a set of matrices $S_\ell\in\mathbb{R}^{m\times n}$ such that 
\[\s = \left\{\sum_{\ell=1}^k\alpha_\ell S_\ell \Bigg| \alpha_1,...,\alpha_k\in\mathbb{R}\right\},\]
we define a \textit{representation matrix} $S\in\mathbb{R}^{m\times nk}$ for the subspace $\s$ as the horizontal concatenation of these matrices, 
\begin{equation} \label{eq:s_basis} S := \begin{bmatrix} S_1 & S_2 & \ldots & S_k \end{bmatrix}.\end{equation}
The following lemma is the tool that will allow us to impose structure on controllers. 

\begin{lem}[Appendix of \cite{ferrante2019design}] 
\label{lem:structure}
Let $\s\subseteq \mathbb{R}^{m\times n}$ be a subspace and let $S$ a representation matrix for $\s$. Define $\Upsilon(S)$ as the following set:
\[ \Upsilon(S) := \{Q \in \R^{n \times n} \mid \exists \Lambda \in \psd^k : S(I_k \otimes Q) = S(\Lambda \otimes I_n)\}.\]
The set $\Upsilon(S)$ is a convex set in $Q$ given $S$, given that is the projection of a subspace in $(Q, \Lambda)$.
Then, the following holds for all matrices $(L, R)$ such that $R$ is invertible: 
% \todo{I don't precisely recall. This is definitely only true for invertible $R$. But wasn't there something like $\Upsilon(S)$ is closed under inverses? }
\[  L \in \s, \ R \in \Upsilon(S) \implies L R^{-1} \in \s.\] 
\end{lem}
\begin{rmk}
    The set $\Upsilon(S)$ is a convex subset of the nonconvex set of matrices $R$ such that $L R^{-1} \in \s$
\end{rmk}
% \begin{rmk}
% If the subspace $\s$ can be represented by matrices $S^1 \neq S^2$ (different bases), the sets $\Upsilon(S^1)$ and $\Upsilon(S^2)$ may be nonequal.
%     % The set $\upsilon(S)$ is dependent on the choice of representation matrix 
% \end{rmk}
% \todo{maybe, if we aren't going to address this problem, it's better to make a less big deal of it, i.e. no remark env., but just a sentence 'note that the choice of representation matrix impacts the set $\Upsilon(S)$}
\begin{exmp}
\label{exmp:sparse}
Continuing from Example \ref{exmp:sparse},  a possible sparsity-derived subspace $\s$ and associated matrix structure $R \in \Upsilon(S)$ for \eqref{eq:sparse} is  \cite[Sec. IV.B]{ferrante2019design} 
    \begin{align}
        \text{sp}(\s) &= \begin{bmatrix}
            1 & 1 & 0 \\ 0 & 1  & 1
        \end{bmatrix},  & 
        R &= \begin{bmatrix}
            R_{11} & R_{12} & 0 \\ 0 & R_{22} & 0 \\ 0 & R_{32} & R_{33}
        \end{bmatrix}. \label{eq:r_upsilon}
    \end{align}
    The structure for $R$ in \eqref{eq:r_upsilon} arises because the representation $\exists \Lambda: \ S(I_k \otimes Q) = S(\Lambda \otimes I_n)$ imposes that $R_{13}, R_{21}, R_{23}, R_{31} = 0$.
    If $R \in \Upsilon(S)$ is invertible, then $L \in \s$ implies $L R^{-1} \in \s$.
\end{exmp}

\begin{rmk}
    The above example included a parametrization of $R$ in \eqref{eq:r_upsilon} that does not rely on finding a value of $\Lambda$. The notation $R \in \Upsilon(S)$ will be used throughout this paper when presenting optimization problems, noting that $\Lambda$ can an optimization variable or eliminated as appropriate.
\end{rmk}

% Example B of \cite{ferrante2019design} involves the following sets for sparse control of an 3-state, 2-input system linear system:
% \begin{align}
%     \s &= \left\{\begin{bmatrix}
%         a_1 & a_2 & 0 \\ 0 & a_3 & a_4
%     \end{bmatrix} \mid a \in \R^4\right\} \\
%     \Upsilon(S) &= \left\{\begin{bmatrix}
%         q_1 & q_2 & 0 \\ 0 & q_3 & 0 \\ 0 & q_4 & q_5
%     \end{bmatrix} \mid q \in \R^5\right\}  \label{eq:upsilon_specific}
% \end{align}
% Restricting $R$ to be diagonal would additionally constraint $q_2, q_4 = 0$ in \eqref{eq:upsilon_specific}.

Using the previous representation matrix formulation, we can derive a structured control design method that is convex in its decision variables.
\begin{lem}\label{lem:struc h2}
   \textit{ Given system \eqref{eq:sys_h2} and a subspace $\s$ with representation matrix $S$ from \eqref{eq:s_basis}. If the following \ac{LMI}s in the variables $P, R, L$ 
    % \todo{Let's be consistent in requiring $\Lambda$ or writing $\Upsilon$.}  
    are simultaneously feasible}
    \begin{subequations}    
    \label{eq:h2_structure}
        \begin{align}
        &\begin{bmatrix}
            P - E E^\top & A R + B L \\
            (A R + B L)^\top & R + R^\top - P 
        \end{bmatrix} \succ 0 \label{eq:h2_aa} \\
        & \begin{bmatrix}
            Q & CR +D L\\ (CR +D L)^\top & R + R^\top - P
        \end{bmatrix} \succ 0 \label{eq:h2_ab}\\
        & P \succ 0 \\
        & \text{Tr}(Q) \leq \gamma^2
     \\
       & P \in \psd^n, \ Q \in \psd^q, \ R \in \Upsilon(S), \ L \in \s,
    \end{align}
    \end{subequations}
    % \todo{We also require $R$ to be invertible!} 
\textit{    then $R$ is nonsingular and  $K =L R^{-1} \in \s$ a $\gamma$-suboptimal $H_2$ controller for \eqref{eq:sys_h2}. }
\end{lem}
Suboptimality of the control design scheme in \eqref{eq:h2_structure} primarily arises through the use of Lemma \ref{lem:structure} as a convex method to design structured controllers.
% Scherer Psatz for Matrices \cite{scherer2006matrix}

%% file: sections/acronym.tex
\begin{acronym}
% \acro{BSA}{Basic Semialgebraic}
\acro{DDC}{Data Driven Control}

% \acro{GAS}{Globally Asymptotically Stable}

% \acro{CSP}{Correlative Sparsity Pattern}

\acro{LMI}{Linear Matrix Inequality}
\acroplural{LMI}[LMIs]{Linear Matrix Inequalities}
\acroindefinite{LMI}{an}{a}

% \acro{LQR}{Linear Quadratic Regulator}
% \acroplural{LMI}[LMIs]{Linear Matrix Inequalities}
% \acroindefinite{LQR}{an}{a}

% \acro{LP}{Linear Program}
% \acroindefinite{LP}{an}{a}
% \acro{OCP}{Optimal Control Problem}

% \acro{ODE}{Ordinary Differential Equation}

% \acro{POP}{Polynomial Optimization Problem}

\acro{PSD}{Positive Semidefinite}

\acro{QMI}{Quadratic Matrix Inequality}

% \acro{PDE}{Partial Differential Equation}

\acro{SDP}{Semidefinite Program}
\acroindefinite{SDP}{an}{a}

\acro{SOS}{Sum of Squares}
\acroindefinite{SOS}{an}{a}

% \acro{WSOS}{Weighted Sum of Squares}

\end{acronym}

%% file: sections/data_driven_structure.tex
\section{Data-Driven Structured Control}
\label{sec:data_driven}

%\urg{Attn Jaap: Please fix the equations and the margins.}
In the previous problem formulated in    \eqref{eq:h2_structure}, we assumed full knowledge of the system matrices. We now consider the case where the system matrices are not known a priori. To be precise, we assume that the matrices $C$, $D$ and $E$ are chosen by the designer as part of the design goal. Hence, we assume that the state and input matrices $A$ and $B$ have to be determined from measurements. 

In order to do so, we will assume that we have access to input and state measurements, and the state is also driven by unobserved bounded process noise.

%This section will define the set of data-consistent plants and the set of plants that satisfy \eqref{eq:h2_ab} with respect to the common law $u=Kx$ (for which $K \in \s$). It will then apply a Matrix S-Lemma to ensure that all data-consistent plants are simultaneously stabilized. The plant matrices   $A, B$ are a-priori unknown (though noise bounds and collected data is available), but it is assumed that $C, D, E$ are all known matrices.

\subsection{Consistent systems}

A sequence of observations for the time window $t = 1, ...,T$ is collected from the true system $(A_\ast, B_\ast)$  with an unknown process noise $w$:
\begin{equation}
x(t+1) = \textstyle A_\ast x(t) + B_\ast u(t) + w(t). \label{eq:stab_W}
    \end{equation}

Such observations could arise from a single trajectory or from multiple independent experiments.
These observations are assembled 
%into a dataset $\dc = \{x(t), u(t), \delta x(t)\}_{t=0}^{T-1}$.
%The data in $\dc$ is collected 
into the following matrices:
% \todo{I like using $U_-$ and $W_-$ to be consistent, but can go either way. }
    
% This data is collected into matrices $(\bx_-, \bu)$ }
% A sampling process acquires a \rev{corrupted} set of observations $(\bx, \bu, \bth)$ \rev{that could have arisen} from a $T$-time-step trajectory of an unknown \ac{LPVA} \eqref{eq:LPVA} system  $(A(\theta), B)$ 
\begin{align}
\label{eq:data}
    \begin{array}{cccccrl}
        \xn & := & [ x(0) & x(1) & \ldots & x(T-1)] & \in \R^{n \times T}  \\
        \bu & := & [ u(0) & u(1) & \ldots & u(T-1)] & \in \R^{m\times T} \\ 
        \bw & := & [ w(0) & w(1) & \ldots & w(T-1)]  & \in \R^{n \times T}\\
        \xp & := & [ x(1) & x(2) & \ldots & x(T) ] & \in \R^{n \times T}. 
    \end{array}
\end{align}
This yields the expression
\begin{equation} \label{eq:data-eq}
    \xp = A_\ast \xn + B_\ast \bu + \bw.
\end{equation}
We assume that we have access to the matrices $\xp$, $\xn$, and $\bu$. In contrast, the matrix $\bw$ is unknown, but we assume it is bounded as follows: 
\begin{asmp}\label{asmp:noise}
The matrix $\bw$ collecting the process noise signal $w(\cdot)$ satisfies 
\begin{equation}
    \begin{bmatrix}
        I \\ \bw^\top
    \end{bmatrix} \begin{bmatrix}
        \Phi_{11} & \Phi_{12} \\
        \Phi_{12}^\top & \Phi_{22}
    \end{bmatrix} \begin{bmatrix}
        I \\ \bw^\top
    \end{bmatrix}^\top \succeq 0. \label{eq:qmi_noise}
\end{equation}
for some $\Phi \in \psd^{n+T}$ with $\Phi_{11} \in \psd_+^n$ and $-\Phi_{22} \in \psd_{++}^T$. \hfill $\Box$
\end{asmp}
The expression in \eqref{eq:qmi_noise} is a \ac{QMI} constraint on $\bw$, which means that $\bw$ is contained in a matrix ellipsoid. 

\begin{rmk}
Noise models from \eqref{eq:qmi_noise} capture the case where $\bw$ has bounded energy, is a confidence interval of a Gaussian distribution, or can be conservatively used to capture element-wise norm bounds \cite{van2023quadratic}. In particular, a recorded trajectory of length $T$ with an per-time process noise bound of $\epsilon$ ($\forall t \in 0..T-1: \ \norm{w(t)}_2 \leq \epsilon)$ can be overapproximated by \eqref{eq:qmi_noise} using a $\Phi$ matrix of 
\begin{align}
    \Phi = \begin{bmatrix}
        T \epsilon I_n & \0 \\ \0 & -I_T
    \end{bmatrix}. \label{eq:noise_bound_eps}
\end{align}
\end{rmk}

% Under Assumption 1, 
% any system $(A,B)$ for which there exists $\bw$ such that  is compatible with the data. 
We now aim to form an expression to bound the consistent plant matrices $(A, B)$ under the bound on $\bw$ from Assumption 1.
Combining the dynamics relation from \eqref{eq:data-eq} and the noise bound in \eqref{eq:qmi_noise}, we can define a matrix $\Psi \in \psd^{2n+m}$ as
\begin{align}
\label{eq:psi_matrix}
    \Psi := \begin{bmatrix}
        I & \xp \\
        \0 & -\xn \\
        \0 & -\bu
    \end{bmatrix}^\top \Phi \begin{bmatrix}
        I & \xp \\
        \0 & -\xn \\
        \0 & -\bu
    \end{bmatrix}.
\end{align}
and denote the set of systems $(A, B)$ that are compatible with the data $\{\xp,\xn,\bu\}$ as
\begin{align}
    \Sigma_\dc = \left \{(A, B) \ \Bigg| \begin{bmatrix}
        I \\ A^\top \\ B^\top
    \end{bmatrix}^{\!\top} \Psi \begin{bmatrix}
        I \\ A^\top \\ B^\top
    \end{bmatrix}\succeq 0  \right\}. \label{eq:data_consistency}
\end{align}
We refer to $\Sigma_{\dc}$ as the set of plants (matrices $(A, B)$) \textit{consistent with the data}. By assumption, this set is nonempty, as the true plant $(A_\ast, B_\ast)$ is a member of $\Sigma_\dc$.
Recall our objective of designing a structured $\gamma$-suboptimal $H_2$ controller for \eqref{eq:sys_h2}. We now reason as follows: Given that we can not distinguish the systems in $\Sigma_\dc$ on the basis of the data $\{\xp,\xn,\bu\}$, we can  guarantee our objective only if $K\in\s$ and $K$ is a $\gamma$-suboptimal $H_2$ controller for all $(A,B)$ compatible with the data. With these elements in place, we therefore aim at resolving whether there exist $P, R, L, \Lambda$ such that \eqref{eq:h2_structure} is feasible for all $(A,B) \in\Sigma_\dc$.
% The formulation in \eqref{eq:data_consistency} can be extended to include quadratically bounded input and state noise (errors-in-variables) using the set description in \cite{bisoffi2024controller}.
% \todo{I don't think this can be done without conservatism. Is this a useful note?  \urg{It is nonconservative up to the typical QMI framework: tight with single quadratic constraint, loose with multiple. }}

\subsection{Controlled set}
Note that, in \eqref{eq:h2_structure}, only \eqref{eq:h2_aa} depends on $(A,B)$. In order to efficiently test whether \eqref{eq:h2_aa} holds for all $(A,B)\in\Sigma_\dc$, we  first reformulate the former as a QMI of the same form as \eqref{eq:data_consistency}.

By the application of a Schur complement \eqref{eq:h2_aa} can be written as
\[
 P  - (AR + BL)(R+R^\top - P)^{-1}(AR + BL)^\top \succ E E^\top.\] 
Rearranging the terms yields
\[ P - E E^\top - \begin{bmatrix}
        A\\ B
    \end{bmatrix}^\top \begin{bmatrix}
        R \\ L 
    \end{bmatrix}(R+R^\top - P)^{-1}\begin{bmatrix}
        R \\ L 
    \end{bmatrix}^\top\begin{bmatrix}
         A \\ B
    \end{bmatrix}\succ 0,  \]
    which can in turn be formulated as 
   \[ \small\begin{bmatrix}
        I \\ A^\top \\ B^\top
    \end{bmatrix}^{\!\!\top}\!\!\begin{bmatrix}
        P- E E^\top\!\! & \0 \\ \0 & \!\!\!\!-\begin{bmatrix}
        R \\ L 
    \end{bmatrix}(R+R^\top- P)^{-1}\begin{bmatrix}
         R \\ L  
    \end{bmatrix}^\top
    \end{bmatrix} \begin{bmatrix}
        I \\ A^\top \\ B^\top
    \end{bmatrix}\succ 0 . \label{eq:sigma_d}
\]

% \urg{check for signs}.
Though this inequality is clearly no longer linear in $R$, $P$, $L$, and $E$, it has the same quadratic structure in $A$ and $B$ as \eqref{eq:data_consistency}. This allows us to employ a matrix-valued S-procedure  \cite[Corollary 4.13]{van2023quadratic} to derive the following: 
\begin{lem} \label{lem:cond} 
    Given data $\{\xp,\xn,\bu\}$ collected under Assumption~\ref{asmp:noise}. Then \eqref{eq:h2_aa} holds for all $(A,B)\in\Sigma_\dc$ if and only if there exist $\alpha\geq 0$, $\beta>0$ such that:
    \[\begin{bmatrix}
        P\!-\! E E^\top \!-\!\beta I\!\! & \0 \\ \0 & \!\!\!\!-\begin{bmatrix}
        R \\ L 
    \end{bmatrix}(R+R^\top- P)^{-1}\begin{bmatrix}
         R \\ L  
    \end{bmatrix}^\top
    \end{bmatrix}-\alpha\Psi\succeq 0. \] 
\end{lem}

All ingredients are now in place to state and prove our main result.
\begin{thm} \label{thm:main}
    Given data $\{\xp,\xn,\bu\}$ collected under Assumption~\ref{asmp:noise} and a subspace $\s$ with representation matrix $S$ from \eqref{eq:s_basis}. Consider the following SDP:
\begin{subequations}    
    \label{eq:h2_structure_dd}
    \begin{align}
         &\inf_{P, R, L, \alpha, \beta, \gamma}  \quad  \gamma  \quad \quad \text{subject to:} \nonumber \\
 & {\small  \begin{pmat}[{..|}]
        P\!-\! E E^\top \!-\!\beta I\!\! & \0  & \0 & \0 \cr
        \0 & \0 & \0 & R \cr
        \0 & \0 & \0 & L  \cr\-
        \0 &   R^\top& L^\top    & R+R^\top \!-\! P \cr
    \end{pmat} }
    \!\!-\! \alpha \begin{pmat}[{|}]
        \Psi & \0 \cr\- \0 & \0 \cr
    \end{pmat}\! \succeq \! 0 \label{eq:pd_data}\\
        & \begin{bmatrix}
            Q & CR +D L\\ (CR +D L)^\top & R + R^\top - P
        \end{bmatrix} \succeq 0 \label{eq:pd_q}\\
        & \text{Tr}(Q) \leq \gamma^2 \\
        & \alpha \geq 0, \beta > 0, \gamma \geq 0   \\
       & P \in \psd^n, \ Q \in \psd^q, \ R \in \Upsilon(S), \ L \in \s.
    \end{align}
    \end{subequations}
    If the program \eqref{eq:h2_structure_dd} is feasible with optimal value $\gamma^*$, 
    then $K = L R^{-1} \in \s$ is a $\gamma^*$-suboptimal $H_2$ controller for all systems $(A,B)\in \Sigma_\dc$.
\end{thm} 
\textbf{Proof:}
Following a Schur complement and Lemma~\ref{lem:cond}, we can see that \eqref{eq:pd_data} is equivalent to the fact that \eqref{eq:h2_aa} holds for all $(A,B)\in\Sigma_\dc$. In turn, we can conclude by Lemma~\ref{lem:struc h2}, that $K = LR^{-1}$ is a $\gamma^*$-suboptimal $H_2$ controller for all systems $(A,B)\in \Sigma_\dc$. 
\hfill$\Box$

Recall that, without making further assumptions on the data or noise, we can not distinguish between the underlying measured system and any other system $(A,B)\in \Sigma_\dc$. Hence, this theorem provides a method to guarantee performance of a  controller for the underlying measured system. Since, in addition the only non-conservative step leading up to Theorem~\ref{thm:main} is the convex relaxation based on Lemma~\ref{lem:structure}, this in turn is also the only source of conservatism in this data-based controller design procedure. 

 %   Note that maximizing the sparsity is a combinatorical problem. Sketch $\ell_1$ regularization. Sketch lasso's, etc. Note that there is a Pareto-front between sparsity and performance. 
%    \item Discuss the nonuniqueness of the basis matrices?

\subsection{Computational Considerations}

In order to implement \eqref{eq:h2_structure_dd} on computational devices, the positive definiteness constraints in  \eqref{eq:pd_data} and \eqref{eq:pd_q} must be replaced by positive semidefinite constraints with respect to a positive tolerance $\eta$ (as in $\succeq \eta I$).

The computational complexity of solving \eqref{eq:h2_structure_dd} can be judged based on the sizes of the (semi)definite constraints in \eqref{eq:pd_data} and \eqref{eq:pd_q}. The matrix in \eqref{eq:pd_data} has size $3n+m$, and the matrix in  \eqref{eq:pd_q} has size $q+n$.

%% file: sections/examples.tex
\section{Numerical Examples}

\label{sec:examples}

MATLAB (2024a) code to generate all examples is publicly available at 
\url{https://www.doi.org/10.3929/ethz-b-000702171}. 
% MATLAB (2021a) code to generate the below examples is publicly available at
% \url{https://github.com/daishuyu/noise-in-observations}. 
Dependencies for these routines include Mosek \cite{mosek92} and YALMIP  \cite{lofberg2004yalmip}. All experiments will employ bounded noise described by the $\Phi$ from \eqref{eq:noise_bound_eps} with appropriate choices of $T$ and $\epsilon$. A numerical tolerance of $\eta = 10^{-3}$ is used to enforce positive definiteness of requisite strict LMI constraints.

% Each example will involve a subspace $\mathcal{S}$ that describes a desired sparsity pattern for the controller $K$.
% \urg{TODO: Coordinated control (some are the same columns) example for Ex. 2}

% Each example will involve a subspace $\mathcal{S}$, and will aim to minimize the stability radius of the closed-loop plant $A+BK$ \urg{is this true? Or should we just do $H_2$? we should write about stability radii in this case.}

\subsection{Example 1: Sparse Control}

This example involves $H_2$-suboptimal control of the following ground-truth plant with $n=3$ states and $m=2$ inputs:
\begin{subequations}
\label{eq:exmp_1}
\begin{align}
    A_\ast &= \begin{bmatrix}
        -0.4095  &  0.4036  & -0.0874 \\
    0.5154  & -0.0815  &  0.1069\\
    1.6715  &  0.7718  & -0.3376
    \end{bmatrix} \nonumber \\
    B_\ast &= \begin{bmatrix}
          0  &      0 \\
   -0.6359 & -0.1098 \\
   -0.0325   & 2.2795
    \end{bmatrix}. \label{eq:exmp_1_sys}
    \intertext{We aim at designing a controller $K$ with the sparsity pattern}
    \text{sp}(\s) &= \begin{bmatrix}
            1 & 1 & 0 \\ 0 & 1  & 1 \end{bmatrix},
\end{align}
in which the corresponding $R$ matrix has a structure observed in \eqref{eq:r_upsilon} from Example \ref{exmp:sparse}.
% \todo{unify the previous with the definition of sparsity pattern.}
The $H_2$-suboptimal control problem is specified by the following matrices:
\begin{align}
    C = &\begin{bmatrix}
        I_3 \\ \0_{2 \times 3}
    \end{bmatrix}, &  D = &\begin{bmatrix}
        \0_{3 \times 2} \\ I_2
    \end{bmatrix},  & E = I_3.
\end{align}
\end{subequations}
In other words, we aim at minimizing the gain between the noise and the signal $\begin{bmatrix} x^\top & u^\top\end{bmatrix}^\top$.

The numerical experiments consider the following design choices for convexly synthesized structured control via LMI \eqref{eq:h2_structure_dd}:
\begin{enumerate}
    \item  A dense matrix $P$ with $P=R$ and $L \in \R^{m \times n}$; 
    \item   A diagonal matrix $P$ with $P=R$ and $L \in \mathcal{S}$; 
    \item A diagonal matrix $R$ and $L \in \mathcal{S}$;
    \item A structured matrix $R \in \Upsilon(\mathcal{S})$ and  $L \in \mathcal{S}$.
\end{enumerate}

Design 1 performs $H_2$ suboptimal control without applying a structural constraint on the gain matrix $K$.
The focus of this work is on Design 4. Design 1 is included to compare against the structure-free control cost. Designs 2 and 3 add more rigid constraints on the $H_2$ program when applying structured control as compared to our proposed Design 4.
Designs 2 and 3 are used as a reference method to compare our contribution in Design 4, since the problem of data-driven structured control has not previously been considered in prior work.

% \todo{Does this mean that we are interested in the gain between  \urg{yes}}
% A trajectory of length $T=20$ is obtained from system \eqref{eq:exmp_1_sys} as corrupted by individual-sample noise of intensity $\epsilon = 0.1$.
% \todo{as mentioned, I think this might be confusing. Can't we do e.g. given sparsity pattern, do experiments with the same length, but $WW^\top \leq \lambda I$, with different levels of $\lambda$? }

Table \ref{tab:exmp_1_eps} reports suboptimal $H_2$ bounds for this system with increasing $\epsilon$ under fixed $T=20$. Table \ref{tab:exmp_1} similarly reports $H_2$ bounds with increasing $T$ under fixed $\epsilon=0.1$. All sampled trajectories begin at the initial point $x(0) = \begin{bmatrix} 1 & 0&0\end{bmatrix}^\top$. 
% as a function of increasing $T$. 
The column heading $(A_\ast, B_\ast)$ refers to solving \ac{LMI} \eqref{eq:h2_design} with respect to the ground truth system under the restrictions described by the design.
In Table \ref{tab:exmp_1_eps}, different trajectories of length $T=20$ are sampled under the different values of $\epsilon$. 
Increasing $\epsilon$ leads to a greater certified closed-loop $H_2$ bound for each implementation. 
The structured control of Design 4 ($R \in \Upsilon(S)$) reduces this $H_2$ bounds as compared to the existing approaches for structured control in Designs 2 and 3.
In Table \ref{tab:exmp_1}, a single underlying trajectory of length $T=20$ is used with noise $\epsilon=0.1$, in which the data set $\mathcal{D}$ is formed by keeping to the first $k$ samples in the given trajectory (heading $T = k$). In the case where $T=5$, the data-driven \ac{LMI} in \eqref{eq:h2_structure_dd} is infeasible for all Designs. 
It is worth noting that increasing the data length $T$ can lead to more conservative solutions: this is due to the approximation of the individual sample noise bound by a single ellipsoid \cite{bisoffi2024controller}. This lack of monotonic decrease of $H_2$ the bound is observed in Table \ref{tab:exmp_1}, in which the $T=20$ bounds are higher than the $T=6$ bounds for Designs 2-4. Even so, Design 4 maintains a lower $H_2$ bound than Designs 2-3.
\begin{table}[h]
    \centering
    \caption{Closed-loop $H_2$ bounds for the system in  \eqref{eq:exmp_1} v.s. $\epsilon$ with $T=20$}
    \begin{tabular}{l|c c c  c}
        Design & $(A_\ast, B_\ast)$ & $\epsilon=0.05$& $\epsilon=0.1$ & $\epsilon=0.15$  \\ \hline 
         1 (Unstructured) &  2.1537&  2.3448  &  3.0939 & 4.5757\\
         2 ($P=R$ diag.)&  3.5658& 4.6619  &  7.4193       & Infeasible \\
         3 ($R$ diag.)&  3.0089& 3.5997  &  4.9506    &9.1999\\
         \textbf{4} ($R \in \Upsilon(S)$) & 2.9794 & 3.5495  &  4.6806 &  8.9710\\
    \end{tabular}
    
    \label{tab:exmp_1_eps}
\end{table}

\begin{table}[h]
    \centering
    \caption{Closed-loop $H_2$ bounds for the system in  \eqref{eq:exmp_1} v.s. $T$ with $\epsilon = 0.1$}
    \begin{tabular}{l|c c c  c}
        Design & $(A_\ast, B_\ast)$ & $T=6$& $T=10$ & $T=20$  \\ \hline 
          1 (Unstructured) &  2.1537& 2.9911 & 2.8156 &  3.0939\\
          2 ($P=R$ diag.) &  3.5658& 6.3386& 7.0963 &     7.4193\\
         3 ($R$ diag.) &  3.0089& 4.5545 & 4.5044 &   4.9506\\
         \textbf{4} ($R \in \Upsilon(S)$) &2.9794  &4.4036 &  4.4323&  4.6806
    \end{tabular}
    
    \label{tab:exmp_1}
\end{table}

The case of (Design 4, $T=6$) leads to the following example solution ($Q$ omitted for space):
\begin{subequations}
\begin{align}
    P &= \begin{bmatrix}
        2.4636  & -0.2394 &  -2.2023\\
   -0.2394  &  1.6204  &  0.5241\\
   -2.2023  &  0.5241  & 11.0167
    \end{bmatrix} \\
    R &= \begin{bmatrix}
        2.7105 &  -0.0045 &0\\
   0& 1.6720 &0\\
    0 &    1.4480  & 12.9819
    \end{bmatrix}\\
    L &= \begin{bmatrix}
        1.4524  &  0.3110  &  0 \\
   0 &   -0.7219  &  2.8893
    \end{bmatrix} \\
    K &= \begin{bmatrix}
         0.5359  &  0.1875   &0\\
   0 & -0.6245   &  0.2226
    \end{bmatrix} \\
    \alpha &= 8.8405, \qquad  \beta = 1.0018\times 10^{-5}.
\end{align}
\end{subequations}

% Using the $R \in \Upsilon(S)$ approach for structured controlin Design 4 reduces the conservatism (as judged by the $H_2$ norm from Design 1) as compared to  Designs 2 and 3. 

\subsection{Example 2: Sparse Sharing Control}

The second example involves a system with $n=12$ states and $m=6$ inputs. The sparsity pattern in this example is
\begin{subequations}
\label{eq:exmp_2}
\begin{align}
    \text{sp}(\s) =\begin{bmatrix}        
        \1_{3 \times 2} & \1_{3 \times 8} & \0_{3 \times 2} \\
        \0_{3 \times 2} & \1_{3 \times 8} & \1_{3 \times 2}
    \end{bmatrix},
\end{align}
and we aim at finding a sharing control policy ($\sum_{j=1}^m u_j(t) = 0$ as enforced by $\1^\top K = \0$) following the sparsity structure in \eqref{eq:exmp_2}. This sharing control can be enforced by adding the constraint $\1^\top L = \0$ to the $H_2$ LMI programs \eqref{eq:h2_design} and  \eqref{eq:h2_structure_dd}.
% enforced through the subspace constraint $\mathcal{S} = \{K \in \R^{6 \times 12} \mid K^\top \1 = \0\}$. 
% Using the notation $e_{ij}$ for a matrix with a 1 at location $(i, j)$ and zeros everywhere else, the coordinated control subspace $\mathcal{S}$ can be described through the $(m-1)n$ basis elements $\{e_{1j} - e_{ij}\}_{i=2..m, j\in 1..n}$.

Similar to before, we aim at minimizing the gain between the noise and the signal $\begin{bmatrix} x^\top & u^\top\end{bmatrix}^\top$, that is, we take
\begin{align}
    C = &\begin{bmatrix}
        I_{12} \\ \0_{6 \times 12}
    \end{bmatrix}, &  D = &\begin{bmatrix}
        \0_{12 \times 6} \\ I_6
    \end{bmatrix},  & E = I_{12}.
\end{align}

% Data is collected under an individual-sample noise bound of $\epsilon = 0.05$. 
Table \ref{tab:exmp_2_eps} lists structured  $H_2$ bounds  under the Designs when $T=100$ and $\epsilon$ increases. Table \ref{tab:exmp_2} similarly reports structured  $H_2$ bounds for varying $T$ under a noise bound of $\epsilon=0.05$. Tables \ref{tab:exmp_2_eps} and \ref{tab:exmp_2} echo the results of Tables \ref{tab:exmp_1_eps} and \ref{tab:exmp_1} respectively. For this larger example, Designs 2 and 3 are always infeasible (over the course of the presented experiments), while Design 4 can return feasible $H_2$ values (outside of $T=100, \epsilon = 0.08$). The lack of monotonic decrease of the $H_2$ norm is again observed for Design 4 in Table \ref{tab:exmp_2}.
\begin{table}[h]
    \centering
    \caption{Closed-loop $H_2$ bounds v.s. $\epsilon$ with $T=100$}
    \begin{tabular}{l|c c c  c}
        Design & $(A_\ast, B_\ast)$ & $\epsilon=0.03$& $\epsilon=0.05$ & $\epsilon=0.08$\\ \hline 
 1 (Sharing) & 10.9309 & 11.7460 & 13.1020 & 19.6526 \\
 2 ($P=R$ diag.)  & Infeasible & Infeasible & Infeasible  & Infeasible \\
3 ($R$ diag.)   & Infeasible & Infeasible & Infeasible  & Infeasible \\
\textbf{4} ($R \in \Upsilon(S)$) &  17.9403  & 23.5366  &38.1294    & Infeasible\\
    \end{tabular}
    \label{tab:exmp_2_eps}
\end{table}

 % 10.9309   11.7460   13.1020   19.6526
 %       NaN       NaN       NaN       NaN
 %       NaN       NaN       NaN       NaN
 %   17.9403   23.5366   38.1294       NaN

\begin{table}[h]
    \centering
    \caption{Closed-loop $H_2$ bounds v.s. $T$ with $\epsilon = 0.05$}
    \begin{tabular}{l|c c c  c}
        Design & $(A_\ast, B_\ast)$ & $T=50$& $T=70$ & $T=100$  \\ \hline 
          1 (Sharing) &    10.9309 &  13.4573  & 12.5807 &  10.6945\\
         3 ($P=R$ diag.)  &  Infeasible & Infeasible& Infeasible& Infeasible \\
         3 ($R$ diag.)  &  Infeasible & Infeasible & Infeasible  & Infeasible \\
         \textbf{4} ($R \in \Upsilon(S)$) & 17.9403 & 40.4927 & 32.8451  &  38.1294\\
    \end{tabular}

    %T=50: 12.4447   16.6050   16.2448   12.9902
    %T =70: 12.8543   15.9222   15.7173   13.6462
    %T = 100: 10.6945, 11.9047, 11.8049, 10.9309
    
    \label{tab:exmp_2}
\end{table}

\end{subequations}

%% file: sections/conclusion.tex
\section{Conclusion}

\label{sec:conclusion}

This paper presents an $H_2$-suboptimal structured control design method for in a data-driven setting. To be precise, we combined a convex relaxation of the structured $H_2$ controller design problem with the informativity approach to data-driven control. This results in  a semidefinite program given solely in terms of measured state and input data under unmeasured process noise. When feasible, this program yields  a structured controller with a suboptimal $H_2$ gain for the all systems compatible with the data. In turn, the controllers guarantees performance when applied to the underlying system.

\subsection*{Future work}
The \ac{QMI}-based structured control methodology has several possible extensions. 
\ac{QMI} descriptions for other aspects of robust control (including $H_\infty$ or mixed sensitivity methods) may be derived \cite{ferrante2020lmi, van2020noisy}. The development of data-based tests for (structured) design for these properties is an area of future developments. Other aspects include data-driven structured static output feedback, decentralized control, and distributed control. Similarly, a number of different model classes, such as Linear Parameter Varying systems \cite{miller2022data, verhoek2024lpvbiquad}, admit a similar approach to data-driven control and hence are amenable to structured control design. 
% As compared to $H_2$ however, the structured control formulation of $H_\infty$-suboptimal control from \cite{ferrante2020lmi} comes at the cost of adding a bisection on an additional free parameter. 
The $H_2$ synthesis task considered in this work may also be adapted to a continuous-time setting, but we note that there exists a practical difficulty of collecting noisy state-derivative observations (as finite differencing is based on dense discrete-time samples).

As a last area of interest, we mention the design of maximally sparse $\gamma$-suboptimal $H_2$ controllers on the basis of data. The literature provides a number of approaches to relax, or otherwise deal with, the combinatorial number of semidefinite programs that arise from this problem. Of particular interest is mapping out quantitative relations between sparsity, performance and the quality of the data. 

% The data-consistency set $\mathcal{D}$ may also be refined using a multi-block S-Lemma (with further conservatism), or even to take into account Errors-in-Variables noise \cite{bisoffi2024controller}. A vital area for further work is to reduce the conservatism of structured control synthesis while ensuring that computation time remains reasonable.                           

%% file: sections/acknowledgements.tex
\section{Acknowledgements}

% \urg{Optional acknowledgements.}
The authors thank Alexandre Seuret and Mihailo Jovanovi\'{c} for discussions about structured control and data-driven control.